%% file: proc_Arxiv.tex
\documentclass[11pt]{amsart}
\usepackage{geometry}                % See geometry.pdf to learn the layout options. There are lots.
\usepackage{epstopdf}

\usepackage{amsmath}
\usepackage{amssymb}
\usepackage{amsthm}
\usepackage{mathrsfs}	
\usepackage{mathtools}
\usepackage{esint}
\usepackage[caption=false]{subfig}

\usepackage{bm}
\usepackage{bbm}

  \usepackage{paralist}
  \usepackage{graphics} %% add this and next lines if pictures should be in esp format
  \usepackage{epsfig} %For pictures: screened artwork should be set up with an 85 or 100 line screen
\usepackage{graphicx}  \usepackage{epstopdf}%This is to transfer .eps figure to .pdf figure; please compile your paper using PDFLeTex or PDFTeXify.
 \usepackage[colorlinks=true]{hyperref}
\newtheorem{theorem}{Theorem}[section]
\newtheorem{corollary}{Corollary}

\newtheorem{lemma}[theorem]{Lemma}

\theoremstyle{definition}
\newtheorem{definition}[theorem]{Definition}
\newtheorem{remark}{Remark}

\newcommand{\e}{\varepsilon}

\newtheorem{property}{Property}

\newcommand{\R}{\mathbb{R}} %% reali 
\newcommand{\N}{\mathbb{N}} 
\newcommand{\TV}{\mathrm{TV}}
\newcommand{\X}{\mathtt{X}}
\newcommand{\conv}{\mathrm{conv}}
\renewcommand{\L}{\mathcal L}

\DeclareMathOperator{\BV}{BV}

\title[Regularity of entropy solutions] %Use the shortened version of the full title
      {Structure and regularity of solutions to 1d scalar conservation laws}

\author[E. Marconi]{}


 \email{elio.marconi@unibas.ch}
 
\thanks{The author is supported by ERC Starting Grant 676675 FLIRT}

\begin{document}
\maketitle

% Enter the first author's name and address:
\centerline{\scshape Elio Marconi$^*$}
\medskip
{\footnotesize
% please put the address of the first author
 \centerline{Universit\"at Basel, Departement Mathematik und Informatik}
   \centerline{Spiegelgasse 1, 4051 Basel, Switzerland}
} % Do not forget to end the {\footnotesize by the sign }

%\medskip
%
%\centerline{\scshape First-name2 last-name2 and First-name3
%last-name3}
%\medskip
%{\footnotesize
% % please put the address of the second  and third author
% \centerline{ First line of the address of the second author}
%   \centerline{Other lines}
%   \centerline{Springfield, MO 65810, USA}
%}
%
\bigskip

% The name of the associate editor will be entered by an editorial staff
% "Communicated by the associate editor name" is not needed for special issue.
% \centerline{(Communicated by the associate editor name)}

%The abstract of your paper
\begin{abstract}
We consider bounded entropy solutions to the scalar conservation law in one space dimension:
\begin{equation*}
u_t+f(u)_x=0.
\end{equation*}
We quantify the regularizing effect of the non linearity of the flux $f$ on the solution $u$ in terms of spaces of functions with bounded generalized variation.
\end{abstract}

\section{Introduction}
We consider the scalar conservation law in one space dimension:
\begin{equation}\label{E_cl}
\begin{cases}
u_t+f(u)_x=0 & \mbox{in }\R^+\times \R, \\
u(0,\cdot)=u_0(\cdot), &
\end{cases}
\end{equation}
where the flux $f\in C^\infty(\R,\R)$ and the function $u:\R^+_t\times \R_x\rightarrow \R$ is the spatial density of the conserved quantity.
We consider bounded entropy solutions: more precisely we require that $u\in C^0([0,+\infty),L^1_{\mathrm{loc}}(\R))\cap L^\infty(\R^+\times \R)$ satisfies \eqref{E_cl} in the sense of distributions and that for every convex entropy $\eta:\R\to \R$ it holds
\begin{equation*}
\eta(u)_t+q(u)_x\le 0
\end{equation*}
in the sense of distributions, where the entropy flux $q$ is defined up to constants by $q'=f'\eta'$. The well-posedness of the Cauchy problem \eqref{E_cl} in the class of bounded entropy solutions with respect to $L^1_{\mathrm{loc}}$ topology is by now classical. A first consequence is the fact that the $\BV$ regularity of $u$ is propagated in time and this implies that we can describe in a satisfactory way the structure of the entropy solution $u$ if $u_0\in \BV(\R)$.

We are interested in the case $u_0\in L^\infty$, which is included in the classical well-posedness result. 
The first result in this direction is the Oleinik one sided Lipschitz estimate: if the flux is uniformly convex ($f''\ge c>0$), then for every $t>0$, the entropy solution $u(t)\in \BV_{\mathrm{loc}}(\R)$ and the following inequality between measures holds:
\begin{equation}\label{E_Oleinik}
D_x u(t) \le \frac{\mathcal L^1}{ct}.
\end{equation}
On the other hand if $f(w)=\lambda w$ is linear the solution is given by
\begin{equation*}
u(t,x)=u_0(x-\lambda t)
\end{equation*}
so $u(t)$ has the same regularity as the initial datum $u_0$.

Between these two extremal cases it is interesting to discuss if some weaker notion of nonlinearity (compared to uniform
convexity) of the flux has some regularizing effect on the entropy solution $u$.
The literature on this problem is large: several results, even in several space dimensions and for more general weak solutions, have been obtained by means of the kinetic formulation of \eqref{E_cl} and averaging lemmas (see \cite{DPLM_averaging_lemma, LPT_kinetic} and the more recent \cite{GL_frac}). 

In order to get quantitative regularity results we need to quantify the nonlinearity of the flux $f$:
\begin{definition}\label{D_degeneracy}
We say that the flux $f$ has \emph{degeneracy} $\bar p\in\N$ if
\begin{enumerate}
\item $\{f''(w)=0\}$ is finite; 
\item for each $w\in \R$ such that $f''(w)=0$ there exists $p\ge 2$ such that $f^{(p+1)}(w)\ne 0$. Let us denote by $p_w$ be the minimal $p\ge 2$ such that $f^{(p+1)}(w)\ne 0$;
\item $\bar p=\max_wp_w$.
\end{enumerate}
If such a $\bar p$ exists we also say that $f$ has \emph{polynomial degeneracy}.
\end{definition}
It was conjectured in \cite{LPT_kinetic} that if the flux $f$ has degeneracy $p\in \N$, then 
for every $\e,t>0$ the entropy solution $u(t) \in W_{\mathrm{loc}}^{s-\e,1}(\R)$, with $s=\frac{1}{p}$. See \cite{Jabin_scalar} for a result in this direction.
However it seems more convenient to express the regularity of the entropy solution in terms of functions with generalized bounded variation: more precisely let $\Phi:[0,ì+\infty)\to[0,ì+\infty)$ be a convex function such that $\Phi(0)=0$, let $v:\R\to \R$ and $I\subset \R$ be 
an interval. We say that $v\in \BV^\Phi(I)$ if
\begin{equation*}
\TV^\Phi v(I):= \sup_{x_1<\ldots < x_n, \newline x_i \in I} \sum_{i=1}^{n-1} \Phi(|v(x_{i+1})-v(x_i)|)<+\infty.
\end{equation*}
See \cite{MO_BVPhi} for an introduction to these spaces.
If $\Phi$ is not degenerate, i.e. $\Phi(h)>0$ for every $h>0$, a function $v\in \BV^\Phi(\R)$ is a regulated function, i.e. for every $\bar x\in \R$ there exist both $\lim_{x\to \bar{x}^-}v(x)$ and $\lim_{x\to \bar{x}^+}v(x)$. This is actually a property that we have for entropy solutions to \eqref{E_cl} if the flux satisfies this minimal nonlinearity assumption: $\{w:f''(w)\ne 0\}$ is dense in $\R$ (see for example \cite{Tartar_notes}). We say in this case that $f$ is \emph{weakly genuinely nonlinear}. Notice that the available fractional Sobolev regularity of the entropy solution does not imply that it is regulated.
An interesting particular case is $\Phi(w)=w^p$, in this case we denote $\BV^\Phi$ with $\BV^{1\over p}$.
We notice that for every $\e>0$ and $p\ge 1$ it holds $\BV^{1\over p}(\R)\subset W^{\frac{1}{p}-\e,p}(\R)$, see \cite{BGJ_fractional}. 

The use of these spaces in this context started in \cite{BGJ_fractional,CJ_BVPhi} to express the regularity of the entropy solution when the flux is strictly (but not necessarily uniformly) convex.

The case of nonconvex fluxes is addressed in the following theorem and it is the final goal of this note.
When not explicitly written we refer to \cite{M_reg} for more details.
\begin{theorem}\label{MT}
Let $f$ be a flux of degeneracy $p$ and let $u$ be the entropy solution of \eqref{E_cl} with $u_0\in L^\infty(\R)$ with compact support.
Then there exists a constant $C>0$, depending on $\mathcal L^1(\mathrm{conv}(\mathrm{supp} u_0))$, $\|u_0\|_\infty$ and $f$, such that for every
$t>0$, it holds
\begin{equation}\label{E_frac}
u(t)\in \BV^{1/p}(\R) \qquad \mbox{and} \qquad \TV^{1/p} u (t) \le C\left(1+\frac{1}{t}\right).
\end{equation}
\end{theorem}

\subsection{Plan of the paper}
In Section \ref{S_lagr} we introduce the main tool of this analysis: an extension to the non smooth setting of the classical method of characteristics called Lagrangian representation. This notion has been developed in different settings: a preliminary version has been introduced in \cite{BM_scalar} for wave-front tracking approximate solutions, in \cite{BM_continuous} it has been adapted to deal with the case of bounded and continuous initial data, then extended to $L^\infty$ initial data in \cite{BM_structure}. Moreover an extension to systems is given in \cite{BM_system}. In this note we only need to give a representation for solutions with piecewise monotone initial data, therefore we follow \cite{M_reg} where a simplified version of the Lagrangian representation is provided.

In Section \ref{S_length} we present the main novelty of \cite{M_reg} and of this presentation. It is an estimate of the oscillation of the entropy solution between two characteristics in terms of their distance and the nonlinearity of the flux.
This estimate plays the role that the Oleinik estimate \eqref{E_Oleinik} plays in the convex case and does not require any nonlinearity assumption on the flux.

Building on this result, the Lagrangian representation and the argument in \cite{Cheng_speed_BV}, we present in Section \ref{S_Cheng} the main steps for proving the $\BV_{\mathrm{loc}}$ regularity of $f'\circ u$ under the assumption of polynomial degeneracy of the flux. In \cite{Cheng_speed_BV} the same problem is considered in the case of one and two inflection points.

Finally in Section \ref{S_frac} we briefly comment about the proof of Theorem \ref{MT}.

\section{Lagrangian representation}\label{S_lagr}
As mentioned in the introduction, the starting point is a precise description of the behavior of the characteristics. In this section we present the notion of Lagrangian representation, which extends the notion of characteristic to the non smooth setting. 
Our strategy is to prove uniform regularity estimates on a dense class of bounded entropy solutions so it is sufficient to consider the case in which $u$ is the entropy solution of \eqref{E_cl} with $u_0$ continuous, bounded and piecewise monotone.
\begin{definition}
We say that $\X:\R^+_t\times \R_y\to \R$ is a \emph{Lagrangian representation} of the entropy solution $u$ if 
\begin{enumerate}
\item $\X$ is Lipschitz continuous with respect to $t$;
\item $\X$ is increasing and continuous with respect to $y$;
\item $\X(0,y)=y$ for every $y\in\R$;
\item for every $t\ge 0$ it holds
\begin{equation}\label{E_repr_form}
u(t,x)=u_0(\X(t)^{-1}(x)),
\end{equation}
for every $x\in \R\setminus N$ with $N$ at most countable.
\end{enumerate}
\end{definition}
\begin{remark}
Requiring \eqref{E_repr_form} for every $t\ge 0$ we implicitly refer to the $L^1$ continuous representative of $u$ in time. Moreover it follows immediately from the monotonicity of $\X$ with respect to $y$ and \eqref{E_repr_form} that if $u_0$ is piecewise monotone then $u(t)$ is piecewise monotone for $t>0$. In order to define pointwise the solution, we consider in this case the lower semicontinuous representative. In any case it is necessary to remove a countable set of points in \eqref{E_repr_form}: these are the points where the preimage $\X(t)^{-1}(x)$ is not a singleton and they are the points where shocks are located.
\end{remark}

The Lagrangian representation enjoys several other properties. First the characteristics travel with the characteristic speed: more precisely for every $y\in \R$ and for $\L^1$-a.e. $t>0$ it holds
\begin{equation}\label{E_char}
\partial_t\X(t,y)=
\begin{cases}
f'(u(t,\X(t,y))) & \mbox{if } u(t) \mbox{ is continuous at }\X(t,y) \\
\frac{f(u(t,\X(t,y)+))-f(u(t,\X(t,y)-))}{u(t,\X(t,y)+)-u(t,\X(t,y)-)} & \mbox{if }u(t)  \mbox{ has a jump at }\X(t,y)
\end{cases}.
\end{equation}
Two other properties are relevant in the following.

\begin{property}\label{Pr_1}
For every $(\bar t,\bar x)\in (0,+\infty)\times \R$ there exists $y\in \R$ such that $\X(\bar t,y)=\bar x$ and  at least one of the following holds:
\begin{enumerate}
\item for every $t\in [0,\bar t]$, 
\begin{equation*}
u(t,\X(t,y)-)\le u_0(y)\le u(t,\X(t,y)+);
\end{equation*}
\item for every $t\in [0,\bar t]$, 
\begin{equation*}
u(t,\X(t,y)+)\le u_0(y)\le u(t,\X(t,y)-).
\end{equation*}
\end{enumerate}
\end{property}
This is a way to formulate in the nonsmooth case the fact that the smooth solutions are constant along characteristics.

In order to state the next property we need to introduce the notion of admissible boundary.
\begin{definition}
Let $T>0$, $w\in \R$ and $\gamma: [0,+\infty)\to \R$ be a Lipschitz curve. Moreover let $u$ be the entropy solution of \eqref{E_cl}
and denote by $\Omega^\pm=\{(t,x)\in [0,T)\times \R: x \gtrless \gamma(t)\}$.
We say that $(\gamma,w)$ is an \emph{admissible boundary} for $u$ up to time $T$ if the restriction of $u$ to $\Omega^-$ is the entropy solution of the initial boundary value problem
\begin{equation*}
\begin{cases}
u_t+f(u)_x=0 & \mbox{in } \Omega^-, \\
u(0,\cdot)=u_0 &  \mbox{in } (-\infty,\gamma(0)), \\
u(t,\gamma(t))=w & \mbox{in }(0,T),
\end{cases}
\end{equation*}
and similary on $\Omega^+$.
\end{definition}

\begin{property}\label{Pr_2}
For every $(\bar t,\bar x)\in (0,+\infty)\times \R$ there exists $y\in \R$ such that $\X(\bar t,y)=\bar x$ and $(\X(\cdot,y),u_0(y))$ is an admissible boundary of $u$ up to time $\bar t$.
\end{property}

A previous extension of the notion of characteristic to the nonsmooth setting is presented in \cite[Chap.10]{Dafermos_book_4}. The characteristic equation \eqref{E_char} implies that for every $y$, the map $t\to \X(t,y)$ is a generalized characteristic in the sense of Dafermos. Therefore the Lagrangian representation $\X$ can be interpreted as a monotone selection of Dafermos generalized characteristics for which \eqref{E_repr_form}, Property \ref{Pr_1} and Property \ref{Pr_2} hold. See also \cite{AGV_Structure} for a similar use in the case of convex fluxes.

\section{Length estimate}\label{S_length}
In this section we present an estimate that relates the distance between two characteristics with the same value and the 
oscillation of the entropy solution between these characteristics.
A relevant feature is the nonlinearity of the flux function $f$ and we quantify it in the following way: given $w_1\le w_2$ we consider twice the $C^0$ distance of $f\llcorner [w_1,w_2]$ from the set of affine functions on $[w_1,w_2]$:
\begin{equation*}
\mathfrak d (w_1,w_2):= \min_{\lambda \in \R} \max_{\{w,w'\}\subset [w_1,w_2]}\left(f(w)-f(w')-\lambda(w-w')\right)
\end{equation*}
\begin{theorem}
Let $u$ be the entropy solution of \eqref{E_cl} with $u_0$ bounded, continuous and piecewise monotone.
Let $t>0$ and $y_l< y_r$ be such that
\begin{enumerate}
\item $u_0(y_l)=u_0(y_r)=:\bar w$;
\item  $\X(\cdot,y_l)$ and $\X(\cdot,y_r)$ enjoy Property \ref{Pr_1} up to time $t$.
\end{enumerate}
Denote by 
\begin{equation*}
s:=\max\{y_r-y_l, \X(t,y_r)-\X(t,y_l)\}
\end{equation*}
and
\begin{equation*}
w_m:=\bar w \wedge \inf_{( \X(t,y_l), \X(t,y_r))} u(t), \qquad w_M:=\bar w \vee \sup_{( \X(t,y_l), \X(t,y_r))} u(t).
\end{equation*}
Then
\begin{equation}\label{E_le}
s\ge \frac{\mathfrak d (w_m,w_M)t}{\|u_0\|_\infty}.
\end{equation}
\end{theorem}

As a corollary we get a first a priori estimate for the entropy solution $u$. Roughly speaking the argument is the following: suppose for simplicity that $u_0$ has compact support. By finite speed of propagation also the solution at time $t$ will have compact support. The estimate \eqref{E_le} tells that each oscillation between two values $a<b$ must occupy a given amount of space, which is strictly greater than 0 if the flux is not affine between $a$ and $b$. But the total amount of space at our disposal is finite so we get an a priori estimate on the number of oscillations between two given values of an entropy solution on a given bounded interval. From this we can immediately recover the compactness in $L^1_\mathrm{loc}$ of the set of equibounded 
entropy solutions if the flux is weakly genuinely nonlinear, which can be obtained for example by a compensated compactness argument (see \cite{Tartar_notes}). This compactness can be made quantitative by means of $\BV^\Phi$ spaces presented in the introduction.

\begin{corollary}\label{C_cor}
Denote by
\begin{equation*}
\mathfrak{N}(h)=\min_{w\in [-\|u_0\|_\infty,\|u_0\|_\infty]}\mathfrak{d}(w,w+h), \qquad
\Psi:=\conv (\mathfrak{N})
\end{equation*}
and for every $\e>0$ set $\Phi^\e(w)=\Psi(\frac{x}{2})x^{\e}$. Then $\forall t>0$
\begin{equation*}
u(t)\in \BV_{\mathrm{loc}}^{\Phi^\e}(\R).
\end{equation*}
\end{corollary}

\begin{remark}
Notice that $\Phi^\e(h)>0$ for every $h>0$ if and only if $f$ is weakly genuinely nonlinear, i.e. $\{w:f''(w)\ne 0\}$ is dense in $\R$.
\end{remark}

In this procedure the length estimate plays the same role as the Oleinik estimate \eqref{E_Oleinik} in order to deduce that the entropy solution $u(t)\in \BV_{\mathrm{loc}}$ for every $t>0$. 
Unfortunately if we specify this last result with $f(u)=u^2$ we get that for every $t>0$ the entropy solution $u(t)\in \BV^{\frac{1}{2}-\e}$ and therefore Corollary \ref{C_cor} is not optimal. 
More in general in the setting of Theorem \ref{MT}, we get $u(t)\in \BV^{\frac{1}{p+1}-\e}$ instead of the expected $u(t)\in \BV^{\frac{1}{p}}$.

\section{BV regularity of $f'\circ u$}\label{S_Cheng}
In this section we discuss the $\BV$ regularity of the velocity field $f'\circ u$. In order to get a positive result we require that the flux function $f$ has polynomial degeneracy (see Definition \ref{D_degeneracy}).

\begin{figure}[htp] 
\centering 
\begin{minipage}[b]{.4\linewidth} 
\centering 
\def\svgwidth{\columnwidth}
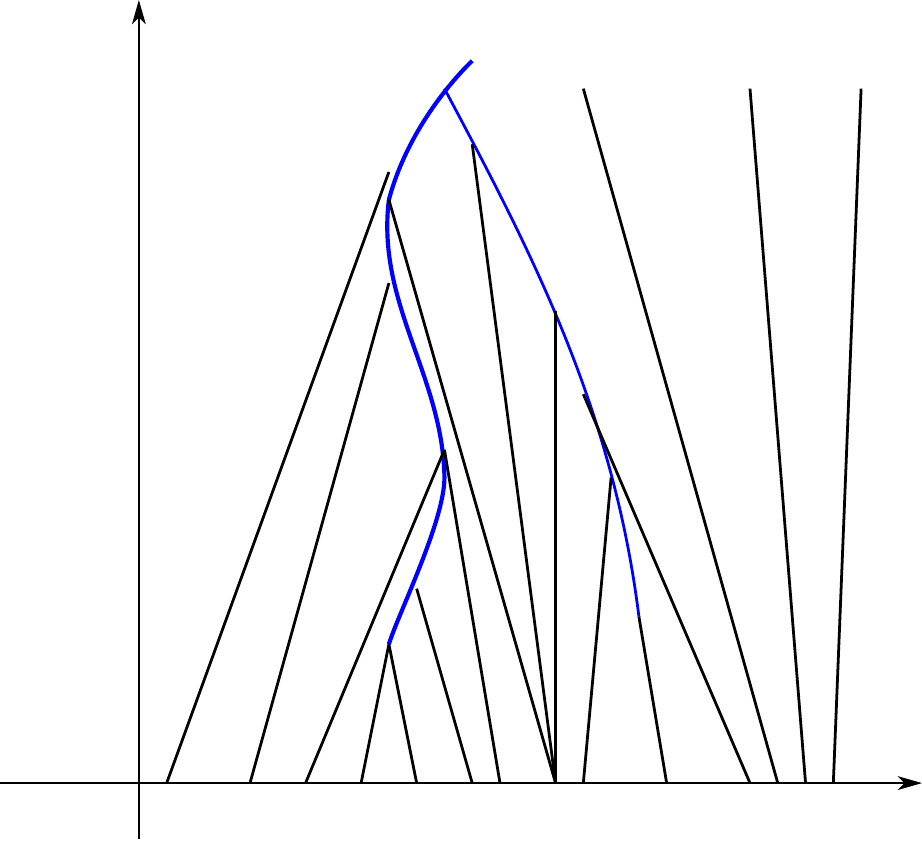
\caption{The characteristics (black) are absorbed by the shocks (blue).}
\label{F_convex}
\end{minipage} 
\hspace{4em}% aumenta o diminuisci secondo come ti sembra meglio 
\begin{minipage}[b]{.45\linewidth} 
\centering 
\def\svgwidth{\columnwidth}
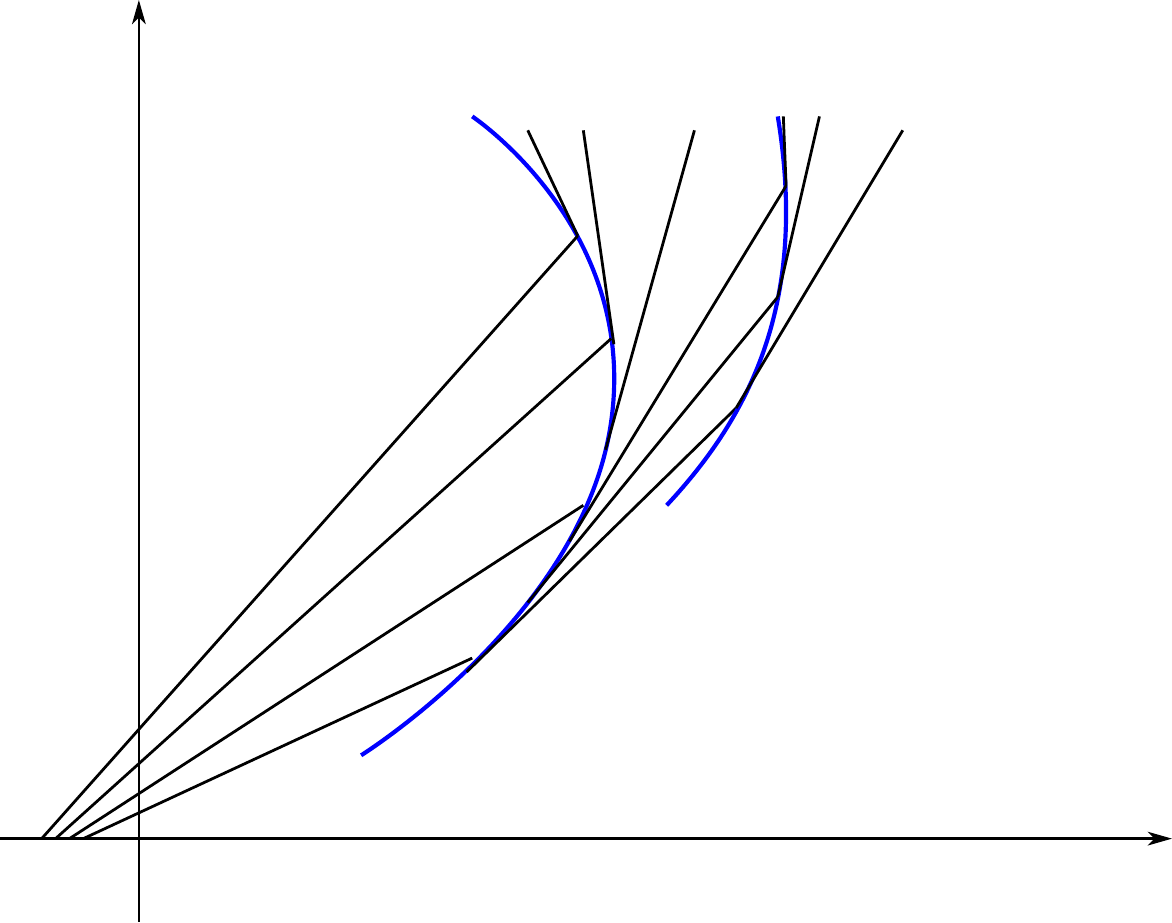
\caption{The characteristics (black) leave the contact discontinuities (blue).}
\label{F_inflection}
\end{minipage} 
\end{figure}

\begin{theorem}\label{T_Cheng}
Let $f$ be as above and $u$ be the entropy solution of \eqref{E_cl} with $u_0\in L^\infty$ and assume that $\mathrm{supp} \,u_0\subset [a,b]$. Then there exists $C$ depending on $b-a, f$ and $\|u_0\|_\infty$ such that for every $t>0$ 
\begin{equation*}
\TV f'\circ u(t) \le C\left(1+\frac{1}{t}\right).
\end{equation*}
\end{theorem}
The details of the proof can be found in \cite{M_reg}. Here we only try to expose the strategy and the role of the tools and the estimates introduced above.
Let us first notice that the situation is much simpler if the flux $f$ is convex. In this case the result follows easily from the structure of the characteristics. The key property is that the characteristics are segments up to the time of the first interaction with other characteristics and two colliding characteristics never split in the future (see Fig. \ref{F_convex}). 
An elementary geometrical constraint and \eqref{E_char} implies that 
\begin{equation}\label{E_speed_O}
D_x f'\circ u(t)\le \frac{\mathcal L^1}{t}
\end{equation}
and the claim easily follows.

This argument does not apply already in the case of fluxes with one inflection point. In this case \eqref{E_speed_O} does not hold and the reason is that two characteristics who interact can split in the future in a contact discontinuity (see Fig. \ref{F_inflection}).
In this case, relying on the precise description of the extremal characteristics in \cite{Dafermos_inflection} and the Lagrangian representation,
the argument in \cite{Cheng_speed_BV} can be made completely rigorous.

The structure of characteristics in the general case is more complicated. It turns out however that it is possible to reduce the general case to the case of fluxes with a single inflection point by means of the length estimate and Property \ref{Pr_2}.
We briefly explain how it can be done:
let $\delta>0$ be the minimal distance between two inflection points of $f$. 
For any $\bar t>0$, thanks to the length estimate \eqref{E_le}, it is possible to find $N\approx C/\bar t$ characteristics 
starting from $y_1<\ldots<y_N$ such that for every $t\in (\bar t/2,\bar t)$ the oscillation of the entropy solution between
two of this characteristics is less than $\delta$. Moreover the constant $C$ depends on the solution only through the 
length of the smallest interval containing the support of $u_0$. 
We therefore obtained $N$ regions in which the range of the solution intersect at most one inflection point of the flux.
The additional difficulty is that in the argument of  \cite{Cheng_speed_BV} we also need to consider the interactions 
of the characteristics with the boundaries of these regions. This can be done interpreting the characteristics as admissible boundaries (Property \ref{Pr_2}) and these are all the ingredients that we need to prove Theorem \ref{T_Cheng}.

\begin{remark}
Actually the $BV$ regularity of $f'\circ u(t)$ can be improved to $SBV$ regularity for every $t\in \R^+\setminus N$ with $N$ countable. See \cite{ADL_SBV} for the case of uniformly convex fluxes, \cite{AGV_SBV} for the extension to the case of strictly convex fluxes and \cite{M_reg} for a proof in the setting of Theorem \ref{T_Cheng}.
\end{remark}

\begin{remark}
The assumption on the flux cannot be removed. In \cite{M_reg} it is provided an example of entropy solution of \eqref{E_cl} in which $f$ has only one inflection point and $f'\circ u$ does not belong to $\BV_{\mathrm{loc}}((0,+\infty)\times \R)$.
\end{remark}

\section{Fractional regularity of the entropy solution}\label{S_frac}
In this last section we deduce Theorem \ref{MT} from Theorem \ref{T_Cheng}. Again, as already noticed in \cite{BGJ_fractional}, the situation is simpler if $f$ is convex. If the flux $f$ has degeneracy $p$, then the inverse function $(f')^{-1}$ is $1\over p$-H\"older and this implies that there exists $C>0$ such that
\begin{equation}\label{E_TVsTV}
\TV^{1\over p} u(t) \le C \TV f'\circ u(t),
\end{equation}
so that Theorem \ref{MT} immediately follows from Theorem \ref{T_Cheng}.

Let us see now how to remove the convexity assumption on $f$: as in the previous section the length estimate allows to  consider only the small oscillations of $u(t)$ and clearly the relevant ones are the oscillations around the inflection points. Therefore it is not restrictive to consider the case $f(u)=u^{p+1}$ with $p$ even.
An estimate like \eqref{E_TVsTV} cannot hold for a generic function $u(t)$ as in the convex case, consider for example a function $v$ which takes only the values $a$ and $-a$ for some $a>0$. In this case $\TV f'\circ v=0$ and $\TV^{1\over p} v$ can be arbitrarily large. This obstruction is excluded taking advantage of the fact that $u(t)$ is the entropy solution of \eqref{E_cl},
roughly speaking if $f'(w_1)\approx f'(w_2)$ the shock between $w_1$ and $w_2$ is not entropic. 
More precisely the following lemma holds.
\begin{lemma}
Let $u$ be the entropy solution of \eqref{E_cl} with $f(u)=u^{p+1}$.
There exists a constant $c>0$ depending on $f$ and $\|u_0\|_\infty$ such that for a.e. $t>0$ and for every $x_1<x_2\in \R$ with $u(t,x_1)\cdot u(t,x_2)<0$ it holds
\begin{equation*}
\TV_{(x_1,x_2)}f'\circ u(t) \ge c|u(t,x_2)-u(t,x_1)|^p.
\end{equation*}
\end{lemma}
By means of this lemma it is not hard to conclude the proof of Theorem \ref{MT}.
\begin{remark}
It has been observed in \cite{CJ_oscillating} that the order $1\over p$ cannot be improved in \eqref{E_frac}.
\end{remark}

% You may incorporate your references as follows in your main tex file.
% Using BibTex is not recommended but can be handled.

\end{document}

%% file: convex.pdf_tex
%% Creator: Inkscape inkscape 0.48.0, www.inkscape.org
%% Accompanies image file 'convex.pdf' (pdf, eps, ps)

\begingroup
  \makeatletter
  \providecommand\color[2][]{%
    \errmessage{(Inkscape) Color is used for the text in Inkscape, but the package 'color.sty' is not loaded}
    \renewcommand\color[2][]{}%
  }
  \providecommand\transparent[1]{%
    \errmessage{(Inkscape) Transparency is used (non-zero) for the text in Inkscape, but the package 'transparent.sty' is not loaded}
    \renewcommand\transparent[1]{}%
  }
  \providecommand\rotatebox[2]{#2}
  \ifx\svgwidth\undefined
    \setlength{\unitlength}{265.53pt}
  \else
    \setlength{\unitlength}{\svgwidth}
  \fi
  \global\let\svgwidth\undefined
  \makeatother
  \begin{picture}(1,0.9191939)%
    \put(0,0){\includegraphics[width=\unitlength]{convex.pdf}}%
    \put(0.09038527,0.88330341){\color[rgb]{0,0,0}\makebox(0,0)[lb]{\smash{$t$}}}%
    \put(0.93398108,0.00545526){\color[rgb]{0,0,0}\makebox(0,0)[lb]{\smash{$x$}}}%
  \end{picture}%
\endgroup

%% file: inflection.pdf_tex
%% Creator: Inkscape inkscape 0.48.0, www.inkscape.org
%% Accompanies image file 'inflection.pdf' (pdf, eps, ps)

\begingroup
  \makeatletter
  \providecommand\color[2][]{%
    \errmessage{(Inkscape) Color is used for the text in Inkscape, but the package 'color.sty' is not loaded}
    \renewcommand\color[2][]{}%
  }
  \providecommand\transparent[1]{%
    \errmessage{(Inkscape) Transparency is used (non-zero) for the text in Inkscape, but the package 'transparent.sty' is not loaded}
    \renewcommand\transparent[1]{}%
  }
  \providecommand\rotatebox[2]{#2}
  \ifx\svgwidth\undefined
    \setlength{\unitlength}{337.53pt}
  \else
    \setlength{\unitlength}{\svgwidth}
  \fi
  \global\let\svgwidth\undefined
  \makeatother
  \begin{picture}(1,0.78671601)%
    \put(0,0){\includegraphics[width=\unitlength]{inflection.pdf}}%
    \put(0.04956887,0.75955285){\color[rgb]{0,0,0}\makebox(0,0)[lb]{\smash{$t$}}}%
    \put(0.94261398,0.02213222){\color[rgb]{0,0,0}\makebox(0,0)[lb]{\smash{$x$}}}%
  \end{picture}%
\endgroup